\documentclass[12pt]{article}
\usepackage[utf8]{inputenc}
\usepackage{amscd,amsthm, amsmath, amsfonts,amssymb}

\newtheorem{theorem}{Theorem}
\newtheorem{corollary}{Corollary}
\newtheorem{lemma}{Lemma}
\theoremstyle{definition}
\newtheorem{definition}{Definition}
\newtheorem{example}{Example}
\newtheorem{remark}{Remark}

\def\F{{\mathbb F}}
\def\K{{\Bbbk}}

\def\Q{{\mathbb Q}}
\newcommand{\baseb}{\hfill $\Box$ \\ } 

\makeatletter

\renewcommand{\@biblabel}[1]{#1.}

\makeatother

\begin{document}


\title{Poincare polynomial for \\ the moduli space   $\overline{{\mathcal M}_{0,n}}({\mathbb C})$  \\ and the number of points in $\overline{{\mathcal M}_{0,n}}({\mathbb F}_q)$} 
\author{N. Amburg$^1$, E. Kreines$^{1,2,3}$, G. Shabat$^{1,4}$}
\date{{\normalsize $^1$ Institute of Theoretical and Experimental Physics, Moscow, Russia \\ $^2$ Lomonosov Moscow State University, Moscow, Russia \\ $^3$ Moscow Institute of Physics and Technology, Dolgoprudny, Russia \\ $^4$ Russian State Humanitarian University, Moscow, Russia}}
 
\maketitle

\begin{abstract}
We obtain a combinatorial proof that the number of points of the space $\overline{{\mathcal M}_{0,n}}({\mathbb F}_q)$ satisfies the recurrent formula for Poincare polynomials of the space $\overline{{\mathcal M}_{0,n}}({\mathbb C})$.

\end{abstract}

\section{Introduction}
Let algebraic \it quasi-projective \rm variety $V$ be defined over the ring 
$\mathbb{Z}$, i.e, it can be defined by a system (possibly, empty) of homogeneous polynomial equations and non-equalities\footnote{i.e., expressions of type  $f\ne0$ (not,  $f>0$)} with integer coefficients. The set  $V(\Bbbk)$ of \it solutions \rm of such system can be considered in the projective space $\mathbf{P}_n(\Bbbk)$ over any field~$\Bbbk$.

Everywhere in this work  ${\mathbb F}_q$ denotes the finite field of $q$ elements. If a complex variety $V(\mathbb{C})$ is \it complete \rm  (i.e. it can be determined without inequalities) and smooth, then its topology is related to the number of its points over finite fields, i.e. with the cardinalities $|V(\mathbb{F}_q)|$ of finite sets  $V(\mathbb{F}_q)$. 
Namely, the (rational) Betti numbers of this variety
$$
\mathrm{b}_k(V):=\dim\mathrm{H}^{j}(V,\mathbb{Q})
$$
can be represented via  \it Hasse-Weil zeta-function $\mathrm{Z}_V$ \rm of the variety $V$. For any prime  $p$ this function is closely related to the generating function for the numbers  $|V(\mathbb{F}_{p^r})|$. More precisely, 
$$\sum_{r=1}^\infty |V(\mathbb{F}_{p^r})|T^r=T\frac{\mathrm{d}}{\mathrm{d}T}
\log\mathrm{Z}_V(T),$$  where
$$
\mathrm{Z}_V(T):=\exp\sum_{r=1}^\infty\frac{|V(\mathbb{F}_{p^r})|T^r}{r}.
$$
In the simplest case  $V=\mathbf{P}_n$, i.e. when  the variety is defined by an  \it empty  \rm set of equations,  we have $|\mathbf{P}_n(\mathbb{F}_q)|=1+q+\ldots + q^n$, and so
$$
\mathrm{Z}_{\mathbf{P}_n}(T):=\exp\sum_{r=1}^\infty\frac{\sum_{k=0}^n p^{kr}T^r}{r}=
\prod_{k=0}^n\exp\sum_{r=1}^\infty\frac{ (p^{k}T)^r}{r}=$$$$
=\prod_{k=0}^n\exp\big(-\log(1-p^kT)\big)=\frac{1}{(1-T)(1-pT)\dots(1-p^nT)},
$$
while the Poincare polynomial of the corresponding complex variety is equal to
$$
\sum_{k=0}^{2n}\mathrm{b}_k(\mathbf{P}_n(\mathbb{C}),\mathbb{Q})t^k=
1+t^2+t^4+\ldots+t^{2n}.
$$
The relations between Hasse-Weil zeta-function and Poincare polynomial is easy to formulate for varieties  \it without odd cohomology\rm, i.e. such that their Betti numbers with odd indexes are zero. Projective spaces and compactified moduli spaces considered in this work are of the above type. In the case when  $V$ has this property, by \cite{Deligne1974} we have
$$
\mathrm{Z}_V=\frac{1}{\prod_{j=1}^{\dim V}P_{2j}},
$$
where the polynomials $P_{2j}\in\mathbb{Z}[T]$ enjoy the following property: its complex roots lie at the circle $|T|=p^{-j}$ and
$
\deg P_{2j}=\mathrm{b}_{2j}(V). 
$

Varieties  $V$ satisfying
$
|V(\mathbb{F}_q)|\in\mathbb{Z}[q],
$
are widely investigated, see, for example,  \cite{BogaartEdixhoven2005}. 

The moduli varieties $\overline{\mathcal{M}_{0,n}}$ which constitute the subject of our considerations are also well known, see, for example, \cite{DM,EtingofCo}. Their cohomology are computed. For example, in \cite{Getzler1995} the following curious relation is established on the base of operada theory: the generating functions 
$$
g(x,t):=x-\sum_{n=2}^\infty\frac{x^n}{n!}\sum_{i=0}^{n-2}(-1)^it^{2(n-i-2)}\dim\mathrm{H}_i(\mathcal{M}_{(0,n+1)})=$$$$=x-\frac{(1+x)^{t^2}-(1+t^2x)}{t^2(t^2-1)}
$$
 and 
 $$
 f(x,t):=x+\sum_{n=2}^\infty\frac{x^n}{n!}\sum_{i=0}^{n-2}t^{2i}
 \dim\mathrm{H}_{2i}(\overline{\mathcal{M}_{(0,n+1)}}),
 $$
are inverse to each other with respect to the \it composition:  \rm  $f(g(x,t),t)\equiv x$.

Our paper is devoted to a direct proof of a recurrence relation for the numbers of points   
$|\overline{\mathcal{M}_{(0,n+1)}}(\mathbb{F}_q)|$, which coincides with the recurrence for Poincare polynomials of compactified moduli spaces $\overline{\mathcal{M}_{(0,n+1)}}(\mathbb{C})$, cf. \it Weil conjectures, \rm
\cite{Deligne1974}. This method of computing the cohomology  of the moduli spaces of curves of positive genus $\mathcal{M}_{g,n}$ and $\overline{\mathcal{M}_{g,n}}$ in cases when the number of points polynomially depend on the order of the field  (note that if $g=1$, then the last condition is satisfed only in the interval $n<11$) was successfully applied earlier, see the papers   \cite{BergstromTommasi2007,Geer2015} and others. 

In the paper \cite{Keel}, 1992, S.~Keel provides a recurrence relations for the Poincare polynomials of  
$\overline{{\mathcal M}_{0,n}}({\mathbb C})$:   
\begin{equation} \label{eq:1}
P_3=1,
\end{equation} 
\begin{equation} \label{eq:2}
P_{n+1}=(1+t^2)P_{n}+\frac{t^2}{2}\sum_{j=2}^{n-2}{n \choose j} P_{j+1} P_{n-j+1} \mbox{ if } n\ge 3.
\end{equation} 

We are going to show that the number of points of the Deligne-Mumford compactification of the  moduli space of genus 0 algebraic curves with $n$ marked and numbered points over a finite field  ${\mathbb F}_q$ satisfies the same conditions (\ref{eq:1}) and~ (\ref{eq:2}).

Our paper is organized as follows: in Section 2 we introduce the moduli space ${{\mathcal M}_{0,n}}({\Bbbk})$ of genus 0 algebraic curves over the field ${\Bbbk}$ and describe its Deligne-Mumford compactification. Section 3 contains the basic information concerning Betti numbers, that can be useful. In Section 4 we provide the point-forgetting map and investigate its properties. Section 5 contains the formulation and proof of the main result.

\bigskip

\section{Deligne-Mumford compactification of  
$\overline{{\mathcal M}_{0,n}}({\Bbbk})$.}

\subsection{ Stable curves} 

In this section we provide basic necessary notions on genus 0 algebraic curves over arbitrary fields, in particular, over finite fields.
\begin{definition}
Stable genus 0 curve with $n$ marked points over a field   ${\Bbbk}$ is a finite union of projective lines  
$\rho=L_1\cup L_2\cup \ldots\cup L_p$,   $L_i\cong {\mathbf P}_1 ({\Bbbk})$, where $n$ different points  $z_1,z_2,\ldots,z_n\in \rho$ are marked in such a way that the following conditions are satisfied:

\begin{itemize}
\item for each  point $z_i$ there is a unique line   $L_j$, satisfying $z_i\in L_j$;
 
\item two different lines have either one or zero common points; 

\item a graph corresponding to $\rho$ is a tree, namely projective lines $L_1 ,$ $L_2,$ $\ldots,L_p$ mark the vertices, two vertices are incident to the same edge if the corresponding lines intersect each other;
   
\item the number of special  (marked or intersection) points in each   $L_j$ is greater than or equal to~3. 
\end{itemize}
\end{definition}
\begin{example}$\mathbf{P}_{1}(\mathbb C)$ is a projective line. Hence, a stable curve over   $\mathbb C$ is a tree of projective lines with marked points on them, see Fig.~1.   
\end{example}

\begin{center}
\begin{picture}(150,180)
\qbezier{(10,70)(70,90)(190,100)}
\qbezier{(-10,20)(30,30)(90,30)}
\qbezier{(-5,50)(30,60)(120,70)}
\qbezier{(0,0)(30,30)(120,170)}
\qbezier{(150,20)(140,130)(180,170)}
\put(5,22){\circle*{4}}
\put(2,10){$1$}
\put(40,28){\circle*{4}}
\put(37,15){$2$}
\put(70,30){\circle*{4}}
\put(67,17){$3$}
\put(15,54){\circle*{4}}
\put(10,41){$4$}
\put(95,67){\circle*{4}}
\put(92,55){$5$}
\put(105,90){\circle*{4}}
\put(102,78){$6$}
\put(148,50){\circle*{4}}
\put(152,48){$7$}
\put(165,147){\circle*{4}}
\put(170,145){$8$}

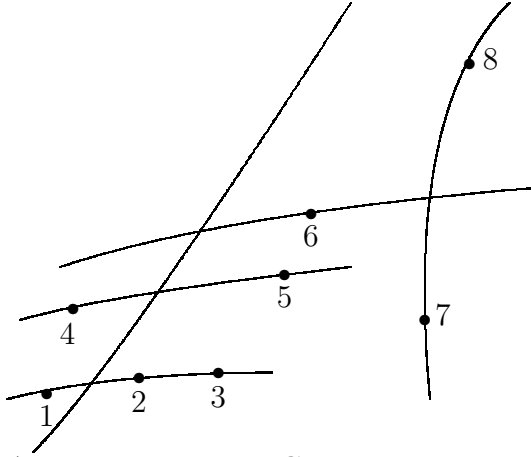
\put(-50,-10){Fig. 1: A stable curve over  $\mathbb C$ with 8 marked points.}
\end{picture}
\end{center}

\bigskip

\begin{remark}
If  $\rho=(L_1,L_2,\ldots,L_p,z_1,z_2,\ldots,z_n)$ is a stable curve, then $n\ge 3$. Thus below we assume   $n\ge 3$, without mentioning this explicitly.
\end{remark}
\begin{example}
Note that $\mathbf{P}_{1}({\mathbb F}_q)={\mathbb F}_q \cup \{\infty\}$ and $\Bigl| \mathbf{P}_{1}({\mathbb F}_q) \Bigr|=q+1$.
Let $(\mathbf{P}_{1}({\mathbb F}_q)\ni z_1,z_2,\ldots,z_n)$ is a stable curve. Then
 $q+1\ge n \ge 3$. 
\end{example}
\begin{example}
Let ${\mathbb F}_2$ be a field of 2 elements,  $n\ge 3$. We consider the union of $n-2$ projective lines $\mathbf{P}_{1}({\mathbb F}_2)$, which are pairwise glued in the $\{0\}$ and $\{\infty\}$, except two which are glued to their neighbors only in one point. The last two constitute the  ``boundary''. We mark the point   $\{1\}$ on each component, and also the remaining points on boundary components. This is a stable curve over  ${\mathbb F}_2$ with $n$ marked points.
\end{example}
\begin{remark}
Note that the inner part of the moduli space   ${\mathcal M}_{0,n} ({\mathbb F}_q)$ is non-empty only in the case  $n\le q+1$, however, for all $n$ and $q$ there exist stable curves with  $n$  
marked points over the field  ${\mathbb F}_q$. 
\end{remark}

\begin{definition}
Let $$\rho=(L_1,L_2,\ldots,L_p,z_1,z_2,\ldots,z_n)\,  \ \ \   \rho'=(L'_1,L'_2,\ldots,L'_p,z'_1,z'_2,\ldots,z'_n)$$ be genus 0 stable curves with  $n$ marked points over an arbitrary field  ${\Bbbk}$. We say that  
$\rho$ and $\rho'$ are equivalent if there exists an isomorphism of algebraic curves   $f:\rho\to \rho'$, such that  $f(z_i)=z'_i$ for all $i=1,\ldots, n$.
\end{definition}
\begin{example} All genus 0 stable curves with $n=3$ marked points over any fixed field   ${\Bbbk}$ are equivalent.
	\end{example}

\subsection{Moduli space  $\overline{{\mathcal M}_{0,n}}({\Bbbk})$}  
\begin{definition}
Let  $n\ge 3$. Deligne-Mumford compactification of the moduli space   ${{\mathcal M}_{0,n}} ({\Bbbk})$ of genus 0 algebraic curves with 
$n$ marked points over a field   ${\Bbbk}$ is the set of equivalence classes of genus 0 stable curves with $n$ marked points over ${\Bbbk}$. We denote this space by $\overline{{\mathcal M}_{0,n}} ({\Bbbk})$. 
\end{definition}
\begin{example} \label{ex2.8}
For any field  ${\Bbbk}$ the moduli space  $\overline{{\mathcal M}_{0,3}}({\Bbbk})$ consists of one point.
\end{example}

\section{Recurrence formula for Betti numbers of  $\overline{{\mathcal M}_{0,n}}({\mathbb C})$}
We provide formal definitions for the notions under consideration.
\begin{definition} Let $X$ be a smooth compact variety. Then the  
$k$-th Betti number  $\beta_k(X)$ is the dimension of the linear space  ${\mathrm H}^k(X; \Q)$, since the cohomology group is a  $\Q$-vector space in this case:
$\beta_k(X)= \dim_{\Q} {\mathrm H}^k(X; \Q)$.
\end{definition}

\begin{definition} Poincare polynomial of  $X$ is the generating function for the Betti number sequence of  $X$:
$$P_X(z)=\beta_0(X)+\beta_1(X)z+\beta_2(X)z^2+\ldots . $$
\end{definition}

\begin{theorem} {\rm \cite{Keel}} \label{keel}
We denote $a_{k}(n)=\beta_{2k}\big(\overline{{\mathcal M}_{0,n}}({\mathbb C})\big)$, i.e., the rank of ${\mathrm H}^{2k}\big(\overline{{\mathcal M}_{0,n}}({\mathbb C}),\mathbb{Q}\big)$. Then the numbers $a_{k}(n)$  satisfy the recurrence relation
\begin{equation} \label{eq:0}
a_{k}(n+1)=a_{k}(n)+a_{k-1}(n)+\frac{1}{2}\sum_{j=2}^{n-2}{n \choose j} \sum_{l=0}^{k-1}a_{l}(j+1)a_{k-1-l}(n-j+1),
\end{equation}
$$
a_{k}(3)= \left\{
\begin{array}{ll}
1, & \mbox{if } k=0; \\
0, & \mbox{if } k\not = 0.
\end{array}
\right.
$$
\end{theorem}

\begin{remark}
It is known (see \cite{Keel}) that $\beta_{2k+1}\big(\overline{{\mathcal M}_{0,n}}({\mathbb C})\big)=0.$
\end{remark}
\begin{corollary}
Using the notation  $P_n(t^2)=\sum_{k=0}^{n-3} a_{k}(n) t^{2k}$, introduced in Section 1  for the Poincare polynomial of   $\overline{{\mathcal M}_{0,n}}({\mathbb C})$, we can rewrite Theorem~{\rm\ref{keel}} in the following form:
$$
P_3=1,
$$
$$
P_{n+1}=(1+t^2)P_{n}+t^2\frac{1}{2}\sum_{j=2}^{n-2}{n \choose j} P_{j+1} P_{n-j+1}, \ \ n>3.
$$ 
\end{corollary}

\section{Point-forgetting map} 

For the completeness we introduce here a classical point forgiving map which we are going to use further to perform the induction. 
\begin{definition}
Point-forgetting map is a natural map 
 $\pi:\overline{{\mathcal M}_{0,n+1}}({\Bbbk})\to \overline{{\mathcal M}_{0,n}}({\Bbbk})$,
cleaning the point  $(n+1)$. If after this cleaning a component   $L_i$ of  $\rho$ contains less than 3 special points, then we change $L_i$ by a point  (in other words, $L_i$ is contracted to a point). 
\end{definition}
\begin{remark}
Note that if   $L_i$ contains less than 3 special points, then, in particular, $L_i$ intersects at most two other components $L_j$ and $L_k$. If $L_i$ intersects exactly two other components, then it is compressed to the intersection point of $L_j$ and $L_k$, and all stability conditions are satisfied. The situation is similar in the other cases. Thus the point-forgetting map is defined correctly. 
\end{remark}
\begin{remark}
Locally, there are three different types of the map $\pi$, they are described in the following example.
\end{remark}
\begin{example}
Different types of point-forgetting map  
 $\pi:\overline{{\mathcal M}_{0,6}}({\Bbbk})\to \overline{{\mathcal M}_{0,5}}({\Bbbk})$ are shown in Fig.~2.
\end{example}

\begin{center}

\begin{picture}(220,100)
\qbezier{(-10,20)(-40,50)(-20,80)}
\put(-20,50){\vector(1,0){50}}
\qbezier{(60,20)(30,50)(50,80)}
\put(-16,27){\circle*{4}}
\put(-25,23){$1$}
\put(-22,37){\circle*{4}}
\put(-31,33){$2$}
\put(-28,47){\circle*{4}}
\put(-37,43){$3$}
\put(-28,57){\circle*{4}}
\put(-37,53){$4$}
\put(-26,67){\circle*{4}}
\put(-35,63){$5$}
\put(-22,77){\circle*{4}}
\put(-30,74){$6$}
 \put(54,27){\circle*{4}}
\put(45,23){$1$}
\put(48,37){\circle*{4}}
\put(39,33){$2$}
\put(42,47){\circle*{4}}
\put(33,43){$3$}
\put(42,57){\circle*{4}}
\put(33,53){$4$}
\put(44,67){\circle*{4}}
\put(35,63){$5$}
\put(127,67){\circle*{4}}
\put(117,65){$6$}
\put(132,37){\circle*{4}}
\put(122,31){$3$}
\put(136,56){\circle*{4}}
\put(134,46){$2$}
\put(150,67){\circle*{4}}
\put(145,71){$4$}
\put(164,80){\circle*{4}}
\put(160,85){$5$}
\put(112,43){\circle*{4}}
\put(109,47){$1$}
\put(218,50){\circle*{4}}
\put(214,39){$3$}
\put(226,56){\circle*{4}}
\put(224,46){$2$}
\put(240,67){\circle*{4}}
\put(235,71){$4$}
\put(254,80){\circle*{4}}
\put(250,85){$5$}
\put(202,43){\circle*{4}}
\put(199,30){$1$}
\qbezier{(140,20)(120,50)(130,80)}
\qbezier{(100,40)(140,50)(170,90)}
\qbezier{(190,40)(230,50)(260,90)}
\put(150,50){\vector(1,0){50}}
\put(0,55){$\pi$}
\put(170,55){$\pi$}

\put(0,0){a)}
\put(170,0){b)}
\end{picture}

\begin{picture}(100,110)
\qbezier{(-65,35)(-15,35)(5,45)}
\qbezier{(-70,20)(-20,30)(10,90)}
\qbezier{(-20,100)(-10,50)(50,70)}
\put(20,55){\vector(1,0){50}}
\put(42,57){$\pi$}
\qbezier{(85,50)(135,50)(155,60)}
\qbezier{(100,80)(110,40)(170,50)}
\put(145,56){\circle*{4}}
\put(141,60){$4$}
\put(102,50){\circle*{4}}
\put(99,42){$3$}
\put(161,48){\circle*{4}}
\put(157,37){$2$}
\put(107,66){\circle*{4}}
\put(111,64){$5$}
\put(142,48){\circle*{4}}
\put(140,37){$1$}
\put(-12,55){\circle*{4}}
\put(-6,52){$6$}
\put(24,64){\circle*{4}}
\put(22,68){$1$}
\put(-5,39){\circle*{4}}
\put(-8,28){$4$}
\put(45,68){\circle*{4}}
\put(40,71){$2$}
\put(-55,34){\circle*{4}}
\put(-59,40){$3$}
\put(-17,92){\circle*{4}}
\put(-23,82){$5$}
\put(50,10){c)}
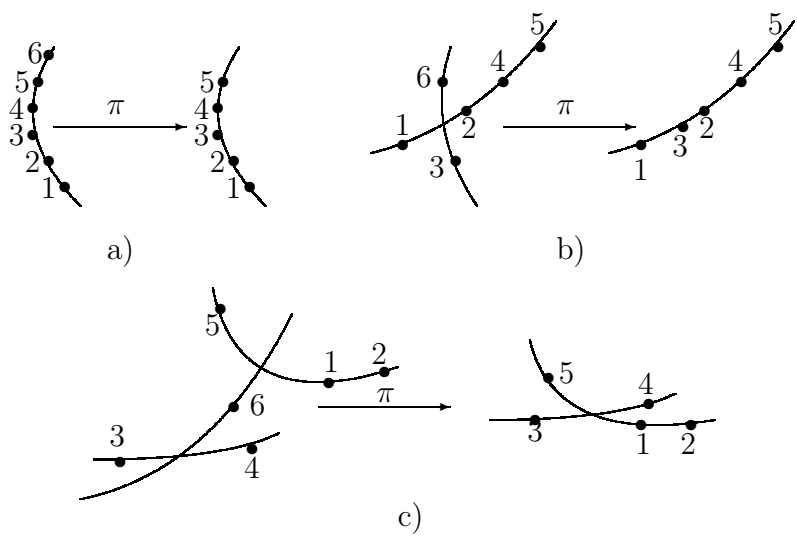
\put(10,-10){Pic. 2. Different types of  
 $\pi:\overline{{\mathcal M}_{0,6}}({\Bbbk})\to \overline{{\mathcal M}_{0,5}}({\Bbbk})$ }
\end{picture}

\end{center}

\bigskip

For a finite field $\K=\F_q$ we compute the number of preimages of any point under the point-forgetting map. For the inner points we have the following.

\begin{lemma} \label{L1}  Let ${\mathbb F}_q$ be a finite field, $|{\mathbb F}_q|=q$,  $\rho\in{{\mathcal M}_{0,{n}}}({\mathbb F}_q)$.
Then 
$$ 
|\pi^{-1}(\rho)|=q+1. 
$$
\end{lemma}

{\bf Proof}.
If $\rho\in{{\mathcal M}_{0,{n}}}({\mathbb F}_q)$, then
$$\rho=(\mathbf{P}_{1}({\mathbb F}_q),z_1,z_2,\ldots,z_n\in \mathbf{P}_{1}({\mathbb F}_q)).$$
Note that $\Bigl| \mathbf{P}_{1}({\mathbb F}_q)\Bigr|=q+1$. The curve $ 
\pi^{-1}(\rho) 
$ can consist of one component. In this case the ``forgotten'' point $z_{n+1}$ could be  any one of   $(q+1-n)$ unmarked points of   $\rho$. This gives  $(q+1-n)$ curves of the space   ${{\mathcal M}_{0,{n+1}}}({\mathbb F}_q)$ (see Fig. 2 {\it a}). Also  ``forgotten'' point $z_{n+1}$ could lie on another component glued in one of   $n$ marked points and contracted by the forgetting map 
(see Fig.  2 {\it b}). This provides another  $n$ curves in ${{\mathcal M}_{0,{n+1}}}({\mathbb F}_q)$. All the obtained curves are different by construction. Thus the full preimage $ 
\pi^{-1}(\rho) 
$
consists of exactly   $(q+1)$ curves.
 \baseb

We denote by $\partial \bar X$ the boundary of   a set   $X$, i.e., a difference between it closure and inner part:   $\partial \bar X=\bar X \setminus \mathring{X}$. For the points in the boundary part the following result holds.
\begin{lemma} \label{L2} 
Suppose that $\rho\in\partial\overline{{\mathcal M}_{0,{n}}}({\mathbb F}_q)$ and the number of irreducible components of  $\rho$ equals $k(\rho)+1$.
Then
$$ 
|\pi^{-1}(\rho)|=(q+1)+qk(\rho). 
$$
\end{lemma}

{\bf Proof}.
Note that if the number of irreducible components of   $\rho$ is $k(\rho)+1$, then the number of their intersection points is smaller by 1, i.e., equals~$k(\rho)$.

As in the proof of the previous lemma the two possibilities occur. The curve $ 
\pi^{-1}(\rho) 
$ has the same number of components as  $\rho$, or not. In the first case the point $z_{n+1}$ could be in any of the  
$$ (k(\rho)+1)(q+1)-n-2k(\rho)$$ unmarked points. The second case also splits: either  $z_{n+1}$ could be in a contracted component glued to one of $n$ marked points (Pic. 2 {\it b}), or  $z_{n+1}$ could be in one of   $k(\rho)$ intersection points (Pic. 2 {\it c}). Together they produce other    $(n+k(\rho))$ curves. Summing the last two numbers we obtain a required result.
\baseb

\section{Recurrence formula for the number of points of   $\overline{{\mathcal M}_{0,n}}({\mathbb F}_q)$} 

Here we prove the main result of the paper. 

\begin{theorem}
\label{Thm}
The number of points
$\Bigl|\overline{{\mathcal M}_{0,n}}({\mathbb F}_q)\Bigr|$ of the space $\overline{{\mathcal M}_{0,n}}({\mathbb F}_q)$ equals  the value of the Poincare polynomial   $P_n(t^2)$
for $\overline{{\mathcal M}_{0,n}}({\mathbb C})$ in 
$t^2=q$:
$$
\Bigl|\overline{{\mathcal M}_{0,n}}({\mathbb F}_q)\Bigr|=P_n(q). 
$$
\end{theorem}

We prove this statement by induction in the number of marked points $n$. We start with some auxhillary lemmas.

\begin{lemma} \label{L3} 
Let $k(\rho)+1$ denote the number of irreducible components of   $\rho\in\overline{{\mathcal M}_{0,{n}}}({\mathbb F}_q)$.
Then
\begin{equation} \label{eq:6}
\Bigl|\overline{{\mathcal M}_{0,{n+1}}}({\mathbb F}_q) \Bigr|=
\Bigl| \overline{{\mathcal M}_{0,n}}({\mathbb F}_q) \Bigr| (q+1)+
\sum_{\rho\in \partial\overline{{\mathcal M}_{0,{n}}}({\mathbb F}_q)}(qk(\rho))
. 
\end{equation}
\end{lemma}

{\bf Proof}. Note that for any curve $\rho\in\overline{{\mathcal M}_{0,{n}}}({\mathbb F}_q)$ the equality   $|\pi^{-1}(\rho)|=
(q+1)+qk(\rho) $  holds. Indeed, if $\rho$ is in the boundary of the moduli space, i.e., it contains several components then the result follows by Lemma \ref{L2}. If  $\rho$ is an inner point, then it consists of just one component, and   $k(\rho)=0$. Then the result follows from Lemma~\ref{L1}. 

Thus
$$ 
\Bigl|\overline{{\mathcal M}_{0,{n+1}}}({\mathbb F}_q)\Bigr|=
\sum_{\rho\in\overline{{\mathcal M}_{0,{n}}}({\mathbb F}_q)} |\pi^{-1}(\rho)|=
\sum_{\rho\in\overline{{\mathcal M}_{0,{n}}}({\mathbb F}_q)}((q+1)+qk(\rho)).$$
Then changing the summands and taking the common multiple  $(q+1)$ we get the equality  $$ \Bigl|\overline{{\mathcal M}_{0,{n+1}}}({\mathbb F}_q)\Bigr| =
\Bigl|\overline{{\mathcal M}_{0,n}}({\mathbb F}_q)\Bigr| (q+1)+
\sum_{\rho\in\overline{{\mathcal M}_{0,{n}}}({\mathbb F}_q)}(qk(\rho))
. 
$$
Since the curves  $\rho\in {\mathcal M}_{0,{n}}({\mathbb F}_q)$ consist of one component, we have    $k(\rho)=0$ for them, and the lemma holds.   
\baseb

For the proof of the next lemma we need some new notations.
\begin{definition} \label{def*}
Let $I=\{i_1,\ldots, i_t\}\subseteq \{1,\ldots,n\}=[n]$ be the set of indexes. 
We denote by  ${\mathcal M}_{0;{I}}({\mathbb F}_q)$ the moduli space of genus 0 algebraic curves over   $\F_q$ with  $k$ marked points numbered by the elements of the set   $I$. By   $\overline{{\mathcal M}_{0;{I}}}({\mathbb F}_q)$ we denote its Deligne-Mumford closure. 
We say that the curve $\rho=(L_1, \ldots, L_p,z_1,\ldots, z_n)\in  \overline{{\mathcal M}_{0,{n}}}({\mathbb F}_q)$  is obtained by gluing $\rho_1$ and $\rho_2$ in the point number  0, if  $\rho_1=(L_{k_1},\ldots, L_{k_s}, z_{i_1},\ldots,z_{i_ t},z_0 )$,   $z_0\in L_{k_l}$ for some $l$; $\rho_2=(L_{k_{s+1}},\ldots, L_{k_p},z_0, z_{j_1},\ldots,z_{j_ {n-t}} )$,    $z_0\in L_{k_r}$ for some $r$,   $k_1,\ldots, k_p$ is a permutation of the set   $\{1,\ldots, p\}$, here $\{j_1,\ldots, j_{n-t}\}=\{1, \ldots, n\}\setminus I$, and the components $L_{k_l}$ and $L_{k_r}$ of $\rho$ intersect.
In this case  $\rho_1 \in \overline{{\mathcal M}_{0;{I}\cup\{0\}}}({\mathbb F}_q)$,  $\rho_2 \in \overline{{\mathcal M}_{0;([n]\setminus {I})\cup\{0\}}}({\mathbb F}_q)$. We denote such gluing by $\rho=\rho_1 *\rho_2$. 
\end{definition}
\begin{remark} Note that for any set $I$ 
$$\overline{{\mathcal M}_{0;{I}}}({\mathbb F}_q) \cong\overline{{\mathcal M}_{0,{|I|}}}({\mathbb F}_q).$$ In particular, 
$|\overline{{\mathcal M}_{0;{I}}}({\mathbb F}_q)| =|\overline{{\mathcal M}_{0,{|I|}}}({\mathbb F}_q)|$, since the choice of the notations for marked points does not change the structure of the moduli spaces.
\end{remark}

\begin{lemma} \label{L4}  
Let $k(\rho)+1$ denote the number of irreducible components of   $\rho\in\overline{{\mathcal M}_{0,{n}}}({\mathbb F}_q)$.
Then
\begin{equation} \label{eq:7}
\sum_{\rho\in\partial\overline{{\mathcal M}_{0,{n}}}({\mathbb F}_q)}k(\rho)=
\frac{1}{2}\sum_{j=2}^{n-2}{n \choose j}\Bigl| \overline{{\mathcal M}_{0,{j+1}}}({\mathbb F}_q)\Bigr|\cdot  \Bigl| \overline{{\mathcal M}_{0,{n-j+1}}}({\mathbb F}_q)\Bigr| 
. 
\end{equation}
\end{lemma}

{\bf Proof}.  Note that in the left hand side of  (\ref{eq:7}) for any curve $\rho\in\partial\overline{{\mathcal M}_{0,{n}}}({\mathbb F}_q)$ we sum the number of intersection points  the irreducible components of this curve.

We consider the expression in the right hand side. Note that any curve  $\rho\in\partial\overline{{\mathcal M}_{0,{n}}}({\mathbb F}_q)$ has at least one intersection point of its components. This point splits the curve into two curves from the spaces with the smaller numbers of marked points.  
Each decomposition   $\{1,\ldots, n\}=\{i_1,\ldots, i_t\}\cup \{i_{t+1},\ldots, i_n\}=I\cup ( [n]\setminus I)$ provides a set of curve pairs   $(\rho_1,\rho_2)$, where $\rho_1\in \overline{{\mathcal M}_{0,I\cup\{0\}}}({\mathbb F}_q)$,  $\rho_2\in \overline{{\mathcal M}_{0,{([n]\setminus I)\cup \{0\}}}}({\mathbb F}_q)$, and marked points are numbered by elements of the sets  $\{i_1,\ldots, i_t,0\}$ and $\{ 0,i_{t+1},\ldots, i_n\}$, correspondingly. By the stability conditions it holds that   $2\le t\le n-2$. Conversely, any pair  of curves  $(\rho_1,\rho_2)$,  $\rho_1\in \overline{{\mathcal M}_{0,I\cup\{0\}}}({\mathbb F}_q)$,  $\rho_2\in \overline{{\mathcal M}_{0,{([n]\setminus I)\cup \{0\}}}}({\mathbb F}_q)$ with marked points numbered by the sets   $\{i_1,\ldots, i_t,0\}$ and  $\{ 0,i_{t+1},\ldots, i_n\}$, respectively, provides one decomposition   $\{1,\ldots, n\}=\{i_1,\ldots, i_t\}\cup \{i_{t+1},\ldots, i_n\}$. Let us consider the set of curves
$${\mathcal X}=\bigcup_{\{i_1,\ldots,i_t\}\subset[n]} \bigcup_{\begin{smallmatrix} \rho_1\in\overline{{\mathcal M}_{0,I\cup\{0\}}}({\mathbb F}_q) \\ \rho_2\in \overline{{\mathcal M}_{0,([n]\setminus I)\cup\{0\}}} ({\mathbb F}_q)\end{smallmatrix}} \rho_1 * \rho_2,$$
where the notation $\rho_1 * \rho_2$ is introduced in Definition~\ref{def*}.
By the above correspondence any curve $\rho$ from $\partial\overline{{\mathcal M}_{0,{n}}}({\mathbb F}_q)$ is included into the set ${\mathcal X}$ exactly $2k(\rho)$ times, since it is counted for any possible split of the curve   $\rho$ into two components and does not depend on the renaming   $\rho_1$ into $\rho_2$ and vice verse. Let us count the number of elements in the set ${\mathcal X}$: 
$$| {\mathcal X}|=\sum_{t=2}^{n-2}  
 {n \choose t}\left|\overline{{\mathcal M}_{0,{t+1}}}({\mathbb F}_q)\right| \cdot \left|\overline{{\mathcal M}_{0,{n-t+1}}}({\mathbb F}_q)\right| .$$
Indeed, for any fixed   $t\in [2,n-2]$ we have the product of the following three numbers: the number of points in the space  $\overline{{\mathcal M}_{0,{t+1}}}({\mathbb F}_q)$, the number of points of the space $\overline{{\mathcal M}_{0,{n-t+1}}}({\mathbb F}_q)$, and the number of choices of $t$ points from $n$. 

Now, taking the common factor   2 in each summand and dividing by it, we obtain on the right-hand side the sum of  $k(\rho)$ over all curves  $\rho$ in  $\partial\overline{{\mathcal M}_{0,{n}}}({\mathbb F}_q)$, i.e., the same amount as in the left-hand side.
\baseb

{\bf Proof of Theorem~\ref{Thm}.}
We prove it by induction in $n$ for each fixed~$q$.

The base of induction is $n=3$. The direct computations and Example~\ref{ex2.8} show that  
$$
\Bigl|\overline{{\mathcal M}_{0,3}}({\mathbb F}_q)\Bigr|=1,
P_3(q)=1 
$$
for all $q$. Thus
$$
\Bigl|\overline{{\mathcal M}_{0,3}}({\mathbb F}_q)\Bigr|=P_3(q), 
$$
for all   $q$, as required.

The induction  step. Assume that for all $m$, satisfying $3\le m \le n$, the statement of Theorem \ref{Thm}  holds, i.e., 
$$
\Bigl|\overline{{\mathcal M}_{0,m}}({\mathbb F}_q)\Bigr|=P_m(q). 
$$
Let us show that it holds for  $m=n+1$.

By Lemma~\ref{L3} we have 
$$
\sum_{\rho\in\partial\overline{{\mathcal M}_{0,{n}}}({\mathbb F}_q)}k(\rho)=
\frac{1}{2}\sum_{j=2}^{n-2}{n \choose j}\Bigl|\overline{{\mathcal M}_{0,{j+1}}}({\mathbb F}_q)\Bigr|\, \Bigl|\overline{{\mathcal M}_{0,{n-j+1}}}({\mathbb F}_q)\Bigr| 
. $$

Taking left hand side of this equality into the formula  (\ref{eq:6}) we obtain that the number of points of the moduli space satisfy the recurrence condition  (\ref{eq:0}). Since initial conditions coincide by the induction base, the number of points of the moduli space satisfy the following conditions:
$$ \left\{ \begin{array}{l}\Bigl|\overline{{\mathcal M}_{0,3}}({\mathbb F}_q)\Bigr|=1,\\
\Bigl|\overline{{\mathcal M}_{0,n+1}}({\mathbb F}_q)\Bigr|=\\ \ \ = (1+q)\Bigl|\overline{{\mathcal M}_{0,n}}({\mathbb F}_q)\Bigr| +q\frac{1}{2}\sum_{j=2}^{n-2}{n \choose j} \Bigl|\overline{{\mathcal M}_{0,j+1}}({\mathbb F}_q)\Bigr|\cdot \Bigl|\overline{{\mathcal M}_{0,n-j+1}}({\mathbb F}_q) \Bigr| , \end{array} \right.
$$
which concludes the proof by Keel theorem  (see Theorem~\ref{keel}). 
\baseb

The work of the first author was done under the partial financial support of RFBR grant 15-01-09242 A, the work of the second author was supported by RFBR grant    15-01-01132. The work of the third author was supported by the Saimons foundation.

\end{document}